\documentclass[12pt]{amsart}
\usepackage{amssymb}

\def\boldn{{\bf n}}

\let\<\langle
\let\>\rangle

\usepackage{euscript}

\newtheorem{theorem}{Theorem}

\newtheorem{lemma}{Lemma}

\newtheorem{slemma}{Star Lemma}

\newtheorem{proposition}{Proposition}

\def\operatorname#1{\mathop{\mathrm{#1}}\nolimits}

\newcommand\C{\mathbb C}
\newcommand\R{\mathbb R}

\def\eqlabel#1{\label{eq#1}}
\def\eqref#1{(\ref{eq#1})}

\def\inf{\operatorname{inf}}

\def\romancap{\operatorname{cap}}
\def\ln{\operatorname{ln}}
\def\vol{\operatorname{vol}}
\def\liminf{\operatorname{liminf}}
\def\mod{\operatorname{mod}}
\def\grad{\operatorname{grad}}

\begin{document}

\title
[Asymptotic Geometry and Conformal Types of C-C Spaces]
{Asymptotic Geometry and Conformal Types of Carnot--Carath\'eodory Spaces}

\author{V.~A.~Zorich}
\address{%
Department of Mathematics and Mechanics\\
Moscow State University\\
119899 Moscow, Russia}


\thanks{Research at MSRI is supported in part by NSF grant DMS-9022140.
The work is supported in part by Russian Basic Research Foundation
project 96-01-00901.}

\begin{abstract}

An intrinsic definition in terms of conformal capacity is proposed for
the conformal type of a Carnot--Carath\'eodory space (parabolic or
hyperbolic). Geometric criteria of conformal type are presented. They
are closely related to the asymptotic geometry of the space at
infinity and expressed in terms of the isoperimetric function and the
growth of the area of geodesic spheres. In particular, it is proved
that a sub-Riemannian manifold admits a conformal change of
metric that makes it into a complete manifold of finite volume if and
only if the manifold is
of conformally parabolic type. Further
applications are discussed, such as the relation between local and
global invertibility properties of quasiconformal immersions (the
global homeomorphism theorem).
\end{abstract}

\maketitle

\section{Carnot--Carath\'eodory Metric and Sub-Riemannian Manifolds}

A Carnot--Carath\'eodory (C-C) metric on a manifold appears when one is
allowed to travel not along arbitrary paths but along {\em
distinguished} ones on the manifold.

The main example of this kind is {\em  polarized Riemannian
manifolds}.

A {\em  polarization} of a manifold $M$ is a subbundle $H$ of the
tangent bundle $TM$.
The polarization $H$ distinguishes a set of directions (tangent
vectors) in $M$ which are usually called {\em  horizontal}.

A piecewise smooth curve in $M$ is called {\em  distinguished}, {\em
admissible}, or {\em  horizontal} with respect to $H$ if the tangent
vectors to this curve are horizontal, i.e., belong to $H$.

The Carnot--Carath\'eodory (C-C) metric $d_{H,\/g}(x_1 , x_2 )$ associated
to the polarization $H$ of the Riemannian manifold $(M, g)$ is defined
as the infimum of the $g$-lengths of the distinguished curves joining
the points $x_1$ and $x_2$.

Thus
$$
d_{H,g} (x_1, x_2) := \inf \,(g\text{-length of } H\text{-horizontal curves
between } x_1 \text{ and } x_2).
$$

This distance obviously satisfies the axioms of a metric, provided
that every two points in $M$ can be connected by a distinguished
curve.

This connectivity property takes place if the Lie
brackets $[[H,H] \ldots]$ of $H$ span the tangent bundle $TM$.

This is the case, e.g., for contact subbundles and for general
non-holonomic (completely non-integrable) distributions.

(About Carnot--Carath\'eodory spaces we refer to \cite{Gro}, where one
can find also a broad bibliography.)

Notice that to define the Carnot--Carath\'eodory metric on $(M,g)$ one
needs a Riemannian structure only on the polarization $H \subset TM$
(but not on the whole tangent bundle $TM$).

A smooth manifold $M$ with a polarization $H \subset TM$ and a
Riemannian structure $g$ on $H$ is refered to as the {\em
sub-Riemannian manifold} $(M, H, g)$.

We assume that the polarization $H$ on $M$ is {\em  generic} in the
sense that $H$ is a smooth distribution of subspaces tangent to $M$,
of equal dimension, and such that Lie brackets (commutators of all
degrees) of $H$ span the tangent bundle $TM$.

We also assume that the polarization $H$ is {\em  equiregular}, i.e.,
that the dimension of the subspace in $T_xM$ generated by the
commutators of fixed degree does not depend on $x \in M$. Then the
tangent bundle $TM$ is filtered by smooth subbundles
\begin{equation}
\eqlabel{1}
H = H_1 \subset H_2 \subset \ldots \subset H_j \subset \ldots \subset
H_k = TM
\end{equation}
such that $H_j$ is spanned by the $j$-th degree commutators of the
fields in $H$.

For example, this is the case if $M$ is a Lie group and $H$ is left
(or right) translation-invariant non-integrable field of tangent
spaces, obtained by translation of the corresponding subspace $t_0M
\subset T_0 M$ tangent to $M$ at the neutral element $0$ of the group.
Here $t_0 M$ is supposed to generate whole the Lie algebra of the
group.

The Heisenberg group $H^n$ with natural (but very unusual)
translation-invariant metric structure is, probably, the most
important example of this kind.

As the simplest example (`flat' in an approrpiate sense), $H^n$ plays
the same role relative to general strictly pseudoconvex
$CR$-manifolds, as $\R^n$ relative to Riemannian manifolds.

(The structure of the Heisenberg group $H^n$ of dimension $n = 2l + 1$
is the same as the
structure of the unit sphere $S^{n}$ in $\C^{l + 1}$. Mostow
observed the close relation between rigidity properties of homogeneous
spaces and quasiconformal mappings of the corresponding sphere at
infinity \cite{Mos}. Complete framework for a theory of quasiconformal
mappings on the Heisenberg group is presented in \cite{KR2}.)

\section{Hausdorff Measure on Carnot--Carath\'eodory Space}

Consider a sub-Riemannian manifold $(M^n, H, g)$ i.e., the smooth
$n$-dimensional manifold $M^n$ with `horizontal' subbundle $H \subset
T M^n$ and a Riemannian metric $g$ on $H$.

We suppose the subbundle $H$ is regular in the sense defined above.
Then $H$ and $g$ induce the {\em  Carnot--Carath\'eodory metric}
$d_{H,g}$ on $M^n$ and turn $M^n$ into {\em  Carnot--Carath\'eodory
manifold} (with horizontal curves as distinguished ones).

Now one can consider the Hausdorff measure of any degree with respect
to this Carnot--Carath\'eodory metric.

It is important to notice that the {\em  metric} (or the {\em
Hausdorff}) {\em  dimension} $m$ of the Carnot--Carath\'eodory space
$(M^n, d_{H,g})$ is usually greater than the topological dimension $n$
of the manifold $M^n$.

Namely,
$$
m = \sum^k_{j=1} j \cdot \text{ rank } (H_j/H_{j-1}),
$$
where $H = H_1 \subset \ldots \subset H_j \subset \ldots \subset H_k =
TM$ is the commutator filtration defined in (1) (we assume $H_0 =
\emptyset$).

The inequality $n < m$ is related to the shape of the
Carnot--Carath\'eodory ball. Carnot--Carath\'eodory metric is highly
non-isotropic and non-homogeneous
in contrast to the usual Euclidean one.

Nevertheless, the Hausdorff measure (volume) induced by
Carnot--Carath\'eodory metric may coincide with the usual Lebesque
measure on a Riemannian manifold or with the invariant Haar measure on
a Lie group.

For instance, the Haar measure on the Heisenberg group $H^n$ coincides
with the Lebesque measure on it $(H^n \approx \R^n)$, as well as with
the Hausdorff measure induced by the above Carnot--Carath\'eodory metric
on $H^n$.

At the same time the metric (Hausdorff) dimension $m$ of the
Heisenberg group (or any domain in it) with respect to the
Carnot--Carath\'eodory metric is equal to $n + 1$.

The corresponding Carnot--Carath\'eodory Hausdorff dimension of a smooth
hypersurface that bounds a domain in $H^n$ is equal to $n$.

A change of the Riemannian structure $g$ on the horizontal bunde $H$
of the sub-Riemannian manifold $(M^n, H,g)$ results in locally
quasiisometric change of the Carnot--Carath\'eodory metric and
thus it does not change the Hausdorff dimension and many other
characteristics of the Carnot--Carath\'eodory space. This is why one
often uses a shorter notation $(M^n, H)$ for the Carnot--Carath\'eodory
structure on the manifold $M^n$.
\section{Horizontal Gradient}

The Riemannian structure $g$ on the horizontal subbundle $H$ is,
nevertheless, very useful for the following definition.

The {\em  horizontal gradient} $\nabla f$ (or $\grad_{H,g} f$)
of a function $f$ on $M^n$ is defined as the unique horizontal vector
such that
$$
\< \nabla f, X \>_g = Xf
$$
for all horizontal vectors $X$.

Here $\<, \>g$ is the inner product with respect to $g$ and $Xf$ is
the Lie derivative of $f$ along $X$ (that is supposed to exist).

The {\em  horizontal normal unit vector} to a hypersurface $\{f = 0\}$
is defined by
$$
\boldn := \frac{1}{\| \nabla f \|} \nabla f.
$$

This is the horizontal normal pointing outward for the domain
$\{ f < 0 \}$.

The vector $ \boldn (x)$ is undefined at points $x \in M^n$ where
$\nabla f = 0$. However, such points form a submanifold of essentially
lower dimension since the horizontal subbundle $H$ is completely
non-integrable. We shall use this remark below without additional comments.

The sub-Riemannian structure $(M^n, H,g)$ on the manifold may be
considered (see e.g., \cite{KR1}) as a limit of the following
Riemannian structures $(M^n, g_\tau)$.

Consider the decomposition $T M^n = H \oplus H'$. Fix $n$ independent
vector  fields $X_1, \ldots, X_d, \ldots X_n$ such that $X_1, \ldots,
X_d$
span $H$ and are orthonormal with respect to the Riemannian
structure $g$ on $H$, and $X_{d + 1}, \ldots , X_n$ span $H'$.

We introduce the Riemannian structure $g_\tau$ on $M^n$ by the condition
that the vectors $\tau X_{d + 1}, \ldots, \tau X_n \ (\tau > 0)$ also
form an orthonormal system.

The one parameter family of Riemannian metrics $g_\tau$ tends to the
singular metric $g$ as $\tau \rightarrow 0$ in the sense that the
length $l_\tau$ of any curve $\gamma \subset M^n$ measured with
respect to $g_\tau$ tends to the Carnot--Carath\'eodory length $l$ of
$\gamma$ with respect to the initial sub-Riemannian structure $(H,g)$
on $M^n$.

This remark allows one to clarify the Carnot--Carath\'eodory counterparts
of some clasical notions and relations that we use below.

For instance, the horizontal gradient is the limit of Riemannian
$g_\tau$-gradients as $\tau \rightarrow 0$.

The classical Fubini-type integral formula %
\begin{equation}
\eqlabel{2}
\int f \, dv = \int_\R \ dt \int_{\{ u = t \} } f \ \frac{d
\sigma}{|\nabla u|}
\end{equation}
(where the domain of integration is foliated by the level surfaces
$\{u = t \}$ of a function $u$) remains valid for sub-Riemannian
manifolds with respect to the induced Hausdorff  volume and area measures
$dv, d\sigma$ and the horizontal grandient $\nabla u$, respectively.

In the special case when $f \equiv 1$ and $u$ is the distance to the
fixed point $O \in M^n$, we  obtain the standard relation
\begin{equation}
\eqlabel{3}
V'(r) = S(r)
\end{equation}
between the volume of the r-ball and the area of its boundary.

The formula \eqref{2} also imply the general relation
\begin{equation}
\eqlabel{4}
dV = \sigma \, dl
\end{equation}
between the volume variation of a regular domain and the area
$\sigma$ of its boundary under the
horizontally-normal variation of the boundary.

The formula \eqref{4} is essentially local: it remains valid for the
local variation of the boundary as well. In this case $\sigma$
should be replaced by the area $d\sigma$
of the variated part of the boundary.

Thus,
\begin{equation}
\eqlabel{5}
d v = d \sigma \ dl
\end{equation}
where $dl$ is the oriented Carnot--Carath\'eodory length of the local
infinitesimal displacement in the direction horizontally normal to the
variated germ of the hypersurface. (One may consider \eqref{5} as the
volume formula for the cylinder.)

Notice that the volume and area integrals in the left- and right-hand
sides of the formula \eqref{2} are meaningful
even when the function under the integration is not
well-defined on some set of volume zero or area zero respectively.

For instance, the Carnot--Carath\'eodory sphere is not a smooth
hypersurface. Thus the horizontal normal
is not well-defined somewhere. But the singularities form a set of area
zero.

\section{Capacity}

In order to characterize conformal types of Riemannian manifolds we
used in \cite{ZK1} some conformal invariants. Namely, the capacity of
an annulus (ring) and the modulus (extremal length) of a family of
manifolds (curves). We need similar special conformal invariants for
the sub-Riemannian case. However, we start with describing the general
notion of capacity (in the spirit of \cite{Maz}) adopted for
sub-Riemannian manifolds.

Consider a smooth $n$-dimensional manifold $M^n$ with a regular
horizontal subbundle $H \subset T M^n$ and a Riemannian metric $g$ on
$H$. Thus, $(M^n, H, g)$ is a sub-Riemannian manifold.

Let $m$ be  its Hausdorff dimension with respect to the induced
Carnot--Carath\'eodory metric, and $dv, d\sigma$
are the volume and area elements, i.e., elements of Hausdorff $m$-measure
for domains and $(m{-}1)$-measure for hypersurfaces in $M^n$.

Let $\Phi(x, \xi)$ be a nonnegative continuous function on the tangent
bundle $T M^n$ which is positive homogeneous of the first degree with
respect to $\xi \in T_x M^n$.

Let $C$ be a compact set in a domain $D \subset  M^n$ and $A = A(C,D)$
be the set $\{ u \}$ of smooth nonnegative functions with compact
support in $D$ such that $u \equiv 1$ in a neighborhood of $C$. We
refer to them as {\em  admissible functions} for the pair $(C,D)$.

The $(p, \Phi)$-{\em capacity} of the set $C$ relative to the domain
$D$ is defined as follows
\begin{equation}
\eqlabel{6}
(p ,\Phi)-\romancap(C,D) := \inf_{u \in A} \int_D \Phi^p (x, \nabla u)
\ dv (x),
\end{equation}
where the infimum is taken over all admissible functions; $\nabla u$
is the horizontal gradient; $dv(x)$ is the volume element of Hausdorff
$m$-measure induced by the Carnot--Carath\'eodory metric of $(M^n, H,g)$.

We need this general notion of capacity only for a special
case of $C,D, p $ and for $\Phi(x, \nabla u) = |\nabla u|$.
Nevertheless, it is reasonable to mention some useful capacity
relations in their general form.

One can rewrite the definition \eqref{6} of the capacity in the
following form for $p > 1$:
\begin{equation}
\eqlabel{7}
(p, \Phi)-\romancap (C,D) = \inf_{u \in A} \left( \int^1_0 \
\frac{dt}{(\int_{\{u = t \}} \Phi^p (x, \nabla u) \frac{d
\sigma}{|\nabla u|})^{1/p-1}}\right)^{1-p},
\end{equation}
(see \cite[section 2.2.2]{Maz}). Here the inner integration is over the
level surface $\{ u = t \}$ of the function $u$ and $d\sigma$ is the
element of area, i.e., of Hausdorff
$(m{-}1)$-measure induced on the level hypersurface by the
Carnot--Carath\'eodory metric.

This new representation \eqref{7} leads to the estimates of capacity
we need below.

Set
\begin{equation}
\eqlabel{8}
P(v) := \inf_{G: |G| \geq v} \int_{\partial G} \Phi(x, \boldn(x))\
d\sigma,
\end{equation}
where the infimum is taken over all domains $G \subset M^n$ with
regular boundary $\partial G$ and such that the volume $|G|$ of $G$ is not
less than $v$; $\boldn (x)$ is the unit vector horizontally normal
to $\partial G$ at the point $x \in \partial G$ directed towards the interior
of $G$.

(The integral in \eqref{8} is well-defined even if $\partial G$ fails to
be smooth on some subset of Hausdorff $(m{-}1)$-measure zero.)

Thus $P(v)$ in \eqref{8} may be considered as a generalized
isoperimetric function. For $\Phi(x, \xi) := |\xi|$ the function
$P(v)$ is the ordinary isoperimetric function of the manifold.
Indeed, in this case for any regular domain $G$ the following relation
is fulfilled
\begin{equation}
\eqlabel{9}
P(v) \leq S,
\end{equation}
where $v$ is the volume (Hausdorff $m$-measure) of $G$ and $S$ is
the area (Hausdorff $(m{-}1)$-measure) of the boundary $\partial G$.

Recall that a function $P: \R_+ \rightarrow \R_+ $ is said to be
the {\em  isoperimetric function} on a manifold (equipped with volume
and area measures) if for every domain $G$ with regular boundary $\partial
G$ the relation \eqref{9} holds.

If the relation \eqref{9} holds for a special family of domains they
say that $P$ is the isoperimetric function for this family.

For example, we need below such an isoperimetric function for the
large Carnot--Carath\'eodory balls that form an exhaustion of a
sub-Riemannian manifold.

At the moment we need function $P$ in \eqref{8}
only for domains $G$ bounded by level surfaces of a function $u$
admissible for the pair $(C,D)$.

By means of this function and of representation \eqref{7} of the
capacity one obtains the following inequality
\begin{equation}
\eqlabel{10}
(p, \Phi)- \romancap (C,D) \geq \left (\int^{|D|}_0 \
\frac{dv}{P^{p/p-1}(v)}\right)^{1-p},
\end{equation}
where $|D|$ is the volume of the domain $D$ and $P$ is defined by
the formula \eqref{8} applied to $G = \{ u \leq t \}$.

The proof of this inequality for a sub-Riemannian manifold is the same
as for the case of a Riemannian one (see \cite[section 2.2.2]{Maz}).

Now turn to the special cases we need.

Let $(M^n, H, g)$ be the sub-Riemannian manifold which is not compact,
but complete with respect to the natural Carnot--Carath\'eodory metric.

Fix a point $0 \in M^n$ that will play the role of
the origin. Let $B(r)$ be the ball of Carnot--Carath\'eodory radius $r$
centered at $0$, $v(r)$ the volume (Hausdorff $m$-measure) of
$B(r)$, $S(r)$ the area (Hausdorff $(m{-}1)$-measure) of the sphere
$\partial B(r)$, $R^b_a := B(b)\backslash \bar{B}(a)$ the annulus
(or ring) domain whose boundary components are spheres of radius $a$
and $b$ $(a < b)$ respectively.

For $\Phi(x, \xi) := |\xi|$ define
\begin{equation}
\eqlabel{11}
\romancap_p R^b_a := (p, \Phi) - \romancap (\bar{B}(a), B(b)).
\end{equation}

For this function $\Phi(x, \xi) = |\xi|$ and for the special choice
$u(x) = (b-a)^{-1}(r(x) - a)$ of the admissible function, where $r(x)$
is the Carnot--Carath\'eodory distance of the point $x$ to the origin
$0$, by means of representation \eqref{7} of the capacity we obtain
the following estimate for the $p$-capacity of the annulus $R^b_a$:
\begin{equation}
\eqlabel{12}
\romancap_p R^b_a \leq \left (\int_a^b \
\frac{dr}{S^{\frac{1}{p-1}}(r)}\right )^{1- p} .
\end{equation}

Along with \eqref{10} applied to the isoperimetric function $P(v)$
of geodesic balls $B(r)$, we finally obtain the following estimates
\begin{equation}
\eqlabel{13}
\left (\int^{v(r_2)}_{v(r_1)} P^{\frac{p}{p - 1}}\right )^{1-p} \leq
\romancap_p R^{r_2}_{r_1} \leq \left (\int_{r_1}^{r_2} S^{-\frac{1}{p
-1}}\right )^{1-p},
\end{equation}
where $v(r_1), v(r_2)$ are volumes of $B(r_1)$ and $B(r_2)$
respectively and $p > 1$ as in \eqref{10} and \eqref{12}.

Remind that in our special case
\begin{equation}
\eqlabel{14}
\romancap_p R^b_a := \inf \int_{M^n} |\nabla u |^{p}(x)\, \, dv(x),
\end{equation}
where the infimum is taken over the smooth functions such that
$u \equiv 0$ in the neighborhood of $B(a)$ and $u \equiv 1$ in the
neighborhood of $M^n \backslash B(b)$, while $dv(x)$ is the element
of the volume (Hausdorff $m$-measure) generated by the
Carnot--Carath\'eodory metric of the sub-Riemannian manifold $(M^n, H, g)$.

Conformal change $\lambda^2 g$ of the Riemannian tensor $g$ on the
bundle $H$ produces local rescaling of the
Carnot--Carath\'eodory lengths element by factor $\lambda$, of horizontal
gradient by factor $\lambda^{-1}$, and of Hausdorff $m$-measure
element by factor $\lambda^m$.

Thus, for $p = m$ the integral in \eqref{14} is invariant with
respect to conformal changes of the sub-Riemannian structure and the
Carnot--Carath\'eodory metric.

We get the following conformally invariant capacity of the annulus
(ring, or condenser) $R^b_a$ in the Carnot--Carath\'eodory space $(M^n, H, g)$:
\begin{equation}
\eqlabel{15}
\romancap_m R^b_a = \inf \int_{M^n} |\nabla u|^m (x) \, d v (x),
\end{equation}
where $m$ is the Hausdorff dimension of the Carnot--Carath\'eodory space
$(M^n, H,g)$, $dv(x)$ is the element of the induced Hausdorff
$m$-measure on $M^n$, $\nabla u$ is the horizontal gradient of the
function $u$, and infimum is taken over all admissible functions $u$
mentioned above.

It is well known (see e.g., \cite{Vai}) that the conformal
capacity of a spherical condenser $R^b_a$ in $\R^n$ (thus $m = n$) is
equal to
\begin{equation}
\eqlabel{16}
\romancap_n R^b_a = \omega_{n-1} \left (\ln \frac{b}{a}\right )^{-n+1},
\end{equation}
where $\omega_{n-1}$ is the area of the unit sphere in $\R^n$.

For the case of the Heisenberg group $H^n$ (instead of $\R^n$),
equipped with the natural $n{-1}$-dimensional polarization $H$ and the
Carnot-Carath\'eodory metric, the conformal capacity $(m = n + 1)$ of
the spherical (with respect to Carnot-Carath\'eodory metric) condenser
$R^b_a$ is equal (see \cite{KR1}) to
\begin{equation}
\eqlabel{17}
\romancap_{n + 1} R^b_a = \omega_{n -1} \left (\ln \frac{b}{a}\right )^{-n},
\end{equation}
where $\omega_{n-1}$ is as above.

After these preparations one can carry over to sub-Riemannian
manifolds the main notions and results described in \cite{ZK1} for
Riemannian manifolds.

Below we formulate the notions and results in full detail, sometimes
with comments, but we omit the proofs parallel to those in \cite{ZK1}.

\section{Ahlfors--Gromov Lemma}

Most of standard geometric characteristics of a space are not
conformally invariant. For instance, the unit ball $B^n \subset  \R^n$
is flat with respect to Euclidean metric, but it admits a conformally
Euclidean metric of constant negative curvature (the Poincar\'e
metric) and at the same time admits a conformal (stereographic)
projection to the sphere $S^n \subset  \R^{n + 1}$ of positive constant
curvature.

Nevertheless, some of geometric relations (that are actually
homogeneous-like or so) preserve or change in a controllable way under
conformal changes of the initial metric.

Recall the following useful Ahlfors--Gromov lemma that we formulated
now for sub-Riemannian manifolds.

\begin{lemma}
If two sub-Riemannian manifolds $(M^n, H, g)$ and $(M^n, H,
\tilde{g})$ are conformally equivalent, then
$$
\int^{\tilde{v}(r_1)}_{\tilde{v}(r_0)} \tilde{P}^{\frac{m}{1-m}} \geq
\int^{r_1}_{r_0} S^{\frac{1}{1-m}}
$$
where $m$ is the common Hausdorff dimension of the Carnot--Carath\'eodory
spaces under consideration; $\tilde{P}$ is any isoperimetric function
of $(M^n, H, \tilde{g})$; $\tilde{v}(r_0), \tilde{v}(r_1)$ are volumes
in $(M^n, H, \tilde{g})$ of any two concentric Carnot--Carath\'eodory
balls $B(r_0), B(r_1)$ of the Carnot--Carath\'eodory space $(M^n, H, g)$;
and $S = S(r)$ is the area (Hausdorff $(m{-}1)$-measure) of the sphere
$\partial B(r)$ in $(M^n, H, g)$.
\end{lemma}

The proof of the lemma for a sub-Riemannian manifold is parallel to
that for the Riemannian case (see e.g., \cite{ZK1}).

The proof shows that Ahlfors--Gromov lemma remains valid even if
$\tilde{P}$ is not a universal $\tilde{g}$-isoperimetric function but
is a $\tilde{g}$-isoperimetric function for $g$-balls only.
\section{Conformal Types of Sub-Riemannian Manifolds}

Consider a noncompact manifold $M^n$ endowed with a horizontal
subbundle $H \subset T M^n $ and a Riemannian metric $g$ on $H$ that
induces the Carnot--Carath\'eodory metric on $M^n$. Let  $m$ be the
Hausdorff dimension of a sub-Riemannian manifold $(M^n, H, g)$ with
respect to the induced Carnot--Carath\'eodory metric.

Let $C$ be a nondegenerate compact set in $M^n$, say, a ball. We are
interested in the conformal capacity $\romancap_m(C, M^n)$.

The manifold $M^n$, or more precisely, the sub-Riemannian manifold
$(M^n, H, g)$, is called {\em  conformally parabolic}  or {\em
conformally hyperbolic} if $\romancap_m (C, M^n) = 0$ or $\romancap_m
(C, M^n) > 0$ respectively.

The capacity relations described in the definition do not depend on
the choice of the compact set $C$ and reflect some conformally
invariant properties of the manifold `at infinity'.

In other words, for the manifold of conformally parabolic type
one has
$$
\lim_{b \rightarrow + \infty} \romancap_m R^b_a = 0,
$$
and for conformally hyperbolic manifold
$$
\lim_{b \rightarrow + \infty} \romancap_m R^b_a > 0
$$
independently of $a > 0$.

The formula \eqref{16} shows that the standard Euclidean space $\R^n$
is of conformally parabolic type.

If we consider the ordinary hyperbolic space modeled on the Euclidean
unit ball, we conclude by means of \eqref{16} that the Lobachevsky
space of constant negative curvature is of conformally hyperbolic
type.

The formula \eqref{17}
shows that the Heisenberg group as a
sub-Riemannian manifold is of conformally parabolic type.

The situation changes if we consider the same group equipped with a
translation-invariant Riemannian structure.

Notice, that the in both these cases the Heisenberg group $H^n$
by itself, endowed with the word-length metric, induces isoperimetric
inequalities with the same isoperimetric function $P(v) = v^{\frac{n +
1}{n}}$.

But the Hausdorff dimension $m$ of the manifold (or rather
Carnot--Carath\'eodory space) $H^n$ is different in the two cases under
consideration.

In the former (sub-Riemannian) case $m = n + 1$ while in the latter
(Riemannian) one $m = n$.

The left-hand-side of \eqref{13} with $p = n$ and $P(v) =
v^{{\frac{n + 1}{n}}}$ shows now that the Heisenberg group $H^n$
equipped with the translation-invariant Riemannian metric is indeed a
manifold of conformally hyperbolic type.

Note that the Carnot--Carath\'eodory metric in many respects is more
adequate for the Heisenberg group than the translation-invariant
Riemannian one. The Carnot--Carath\'eodory metric on $H^n$ is also
translation-invariant but, in addition, this metric follows the
similitudes of the Heisenberg group. It changes by a factor under such
homogeneous transformations.

These two metrics on the Heisenberg group are not conformal or
quasiconformal equivalent even locally.

By the way, every metric space that admits self similitudes (such as
$\R^n$ with Euclidean metric or $H^n$ with Carnot--Carath\'eodory metric)
must be of conformally parabolic type.

\section{Asymptotic Geometry and Conformal Types of Sub-Riemannian
Manifolds}

As it was mentioned above the conformal type of a manifold depends on
the behavior of the manifold at infinity.

We present now an explicit geometric version of this general claim.

We say that a certain property or relation is {\em  realizable in a
class of metrics} if it is valid for some metric of the class.

For instance, any sub-Riemannian manifold $(M^n, H, g)$ can be
realized as a complete one in the conformal class
of its metric.

Below under conformal change of the metric on sub-Riemannian manifold
$(M^n, H, g)$ we mean a conformal change $\tilde{g} = \lambda^2 g$ of
the Riemannian structure on the horizontal bundle $H$. It results in
the conformal change of the corresponding Carnot--Carath\'eodory metric.

Now we are in a position to formulate the following theorem.

\begin{theorem}
Let $(M^n, H, g)$ be a noncompact sub-Riemannian manifold of the
Hausdorff dimension $m$ with respect to the induced
Carnot--Carath\'eodory metric.

The manifold is of conformally parabolic type if and only if any of
the following equivalent conditions is realizable in the class of
complete metrics conformally equivalent to the initial one
\begin{quote}
\begin{itemize}
\item[(i)] $\vol_m (M^n) < \infty$ ,
\item[(ii)] $\int^\infty \ S^{\frac{1}{1-m}} (r) \ dr = \infty$,
\item[(iii)] $\int^\infty \ (\frac{r}{v(r)})^{\frac{1}{m{-}1}} \ dr =
\infty$,
\item[(iv)] $\liminf_{r \rightarrow \infty} \ \frac{v(r)}{r^m} <
\infty$.
\end{itemize}
\end{quote}
\end{theorem}

Here, as above, $v(r)$ and $S(r)$ are the volume and area (i.e., the
Hausdorff $m$-measure and $(m{-}1)$-measure) respectively of the
Carnot--Carath\'eodory ball $B(r)$ and its boundary sphere $\partial B(r)$;
$\vol_m (M^n)$ is the Hausdorff $m$-measure of the whole manifold $M^n$.

The proof is similar to that for the Riemannian manifolds and we omit
it (see \cite{ZK1}).

\section{Canonical Forms of Isoperimetric Function}

In this section we supplement two statements related to possible
variations of the isoperimetric function under conformal changes of
the Carnot--Carath\'eodory metric.

Consider a $g$-spherical exhaustion of a sub-Riemannian manifold
$(M^n, H, g)$, i.e., a system of Carnot--Carath\'eodory balls $B(r)$ of
varying radius $r$ and fixed center $0 \in M^n$. The system is
invariant under the {\em  spherically conformal change} $\tilde{g} =
\lambda^2 (r)g$ of the Riemannian metric $g$ on $H$. Let $m$ be the
Hausdorff dimension of $(M^n, H, g)$ with respect to the induced
Carnot--Carath\'eodory metric.

We denote by $v(r)$, $S(r)$ and $\tilde{v}(r)$, $\tilde{S}(r)$ the
volume (Hausdorff $m$-measure) of the ball $B(r)$ and the area
(Hausdorff $(m{-}1)$-measure) of the sphere $\partial B(r)$ considered in
spaces $(M^n, H, g)$ and $(M^n, H, \tilde{g})$ respectively.
Let $P$ be a nonnegative function on a sub-Riemannian manifold $(M^n,
H, g)$. Consider the class of metrics spherically conformally
equivalent to the given one along with the corresponding spherical
exhaustion of the manifold.

\begin{proposition}
\label{1}
The function $P$ is an isoperimetric function and, moreover, the
maximal one (i.e., $P(\tilde{v}) = \tilde{S}$) for sufficiently large
balls of the spherical exhaustion in a certain metric from this class
if and only if the integrals
$$
\text{ a) } \int^{\bullet} P^{\frac{m}{1-m}} \quad, \quad \text{ b) }
\int^{\bullet} S^{\frac{1}{1-m}}
$$
converge or diverge simultaneously at the upper bound of the domains
of integrands.

Moreover, the new Carnot--Carath\'eodory metric
(induced by $\tilde{g} = \lambda^2 (r) g$) is complete if and only if
$\int^{\bullet} P^{-1} = \infty$.
\end{proposition}

\begin{proposition}
\label{2}
Let $(M^n, H, g)$ be a sub-Riemannian manifold of the Hausdorff
dimension $m$ with respect to the Carnot--Carath\'eodory metric on
$M^n$ induced by the Riemannian structure $g$ on the horizontal
subbundle $H \subset T M^n$.

Consider the class of metrics conformally
equivalent to the given one on $(M^n, H, g)$ along with the geodesic
exhaustion of $M^n$, i.e., the exhaustion by balls with respect to the
corresponding Carnot--Carath\'eodory metric.

Then in this class of conformal metrics and
of the corresponding geodesic exhaustion the maximal isoperimetric
function (for large balls) can be reduced to the following canonical
forms
$$
P(\tilde{v}) = \tilde{v}^{\frac{m-1}{m}} \quad \text{ or } \quad
P(\tilde{v}) = \tilde{v}
$$
according to whether the sub-Riemannian manifold $(M^n, H, g)$ is of
conformally parabolic or conformally hyperbolic type respectively.
\end{proposition}

The Proposition \ref{2} is in a sense incomplete. Roughly speaking,
any function $P$ may occur on the conformally parabolic manifold as
the maximal isoperimetric function for large balls of the geodesic
exhaustion corresponding to a certain metric conformally equivalent to
the initial one.
Besides, it is sometimes possible to get the canonical form of the
isoperimetric function even by the spherically conformal change of the
initial metric.

At last, there is no indication to what extent the Proposition \ref{2}
is invertible.

We omit proofs of Proposition \ref{1} and Proposition \ref{2}. They
are similar to the ones for Riemannian manifolds (cf. \cite{ZK1}).
Notice, however, that Proposition \ref{2} can be essentially developed
in spirit of the paper \cite{ZK2}.

\section{Concluding Remarks}

Asymptotic geometry or, in other words, behavior of the manifold
(space) at infinity sometimes plays a decisive role in global
problems. For instance, it is responsible for the existence of special
solutions of operators
on the manifold. In the case of the Laplace operator the problem is
often closely related to (and for two dimensional manifolds coincides
with) the problem of the manifold type discussed above.

Conformal type of the manifold arises as we consider global
transformations of the manifold or the questions related to the
mappings of one manifold into another.

Liouville-type theorems are examples of this kind. In the classical
form the Liouville theorem claims that there is no bounded
non-constant entire function. In other words, there is no holomorphic
mapping of the plane into a disk. This phenomenon of nonexistence is
of rather general nature and it holds for a much broader class of
mappings. It is related to different conformal types of the Euclidean
plane and the disk (or the hyperbolic space). The former is
conformally
parabolic, while the latter is conformally hyperbolic.

We complete the paper with one more example.

The global homeomorphism theorem (GHT) is the following specifically
multidimensional phenomenon: {\em any locally invertible
quasiconformal mapping} $f: \R^n
\rightarrow \R^n$ {\em is globally invertible provided} $n \geq 3$.

The theorem essentially remains valid for mapings $f: M^n \rightarrow
N^n$ of Riemannian manifolds $(M^n, g_M)$, $(N^n, g_N)$ provided $n
\geq 3$, $\pi_1 (N^n) = 0$, and $(M^n, g_M)$ is of conformally
parabolic type.

There are serious evidence to expect the validity of the GHT for
sub-Riemannian manifolds as well, and with the same condition of
conformal parabolicity on $(M^n, H, g)$.

The initial idea of the proof \cite{Zo1}, which was used in further
generalizations and developments of the GHT
(see \cite{Zo2}), seems to be working for the sub-Riemannian case.
Indeed, the topological part of the proof holds for sub-Riemannian
manifolds as well. Thus, it remains to prove some capacity estimates
equivalent to lemmas 1 and 2 of \cite{Zo1}.

For instance, it is sufficient to prove the following.
\begin{slemma}
Let $D$ be a domain in the Carnot--Carath\'eodory space $(M^n, H, g)$
star-like relative to a point $0 \in D$, i.e., $D$ is formed by rays
or by their intervals originating from $0$ and horizontally normal to
Carnot--Carath\'eodory spheres centered at $0$. The finite intervals end
at points of the part $\triangle$ of the boundary $\partial D$ `visible'
from $0$.

Let $\Gamma_\epsilon $ be the family of all curves in $D$ which join
$\triangle $ with $(\epsilon > 0)$-neighborhood $B(\epsilon)$ of $0$.
If $\mod_m \Gamma_\epsilon = 0$ then
\begin{itemize}
\item[(i)] projection of $\triangle$ along these rays to the sphere
$\partial B(\epsilon) $ does not contain connected components different
from a point;
\item[(ii)] $\partial D$ does not locally divide the manifold $M^n$
provided $n \geq 3$.
\end{itemize}
\end{slemma}

Here $m$ is the Hausdorff dimension of the sub-Riemannian manifold
$(M^n, H, g)$ considered.

For our purposes we may suppose the manifold to be of conformally
parabolic type.

Conformal modulus $\mod_m \Gamma$ of the family $\Gamma$ of curves in
$(M^n, H, g) $ is defined as follows
$$
\mod_m \Gamma := \inf \int_{M^n} \rho^m (x) \ d v(x),
$$
where the infimum is taken over all Borel measurable nonnegative functions
such that $\int_\gamma \rho \geq 1$ for each curve $\gamma \in \Gamma$.

One readily shows (cf. the arguments above
concerning the conformal capacity) that $\mod_m \Gamma$ is a conformal
invariant.

It is closely related to the conformal capacity. For instance,
$$
\romancap_m R^b_a = \mod_m \Gamma^b_a,
$$
where $\Gamma^b_a$ is the family of all curves that  join boundary
components of the spherical condenser $R^b_a$.

Notice that the Star lemma formulated above for the sub-Riemannian
case is true for Riemannian manifolds. In the latter case it can be
proved by means of the slightly modified method that was used in the
proof of lemmas 1 and 2 in \cite{Zo1}.

\subsection*{Acknowledgements.} I wish to thank M.~Gromov, G.~Margulis, and
A.~Schwarz for invitations, fruitful discussions and hospitality during
my visit to University of Maryland, Yale University and University of
California in Spring 1996. I am grateful to B.~Khesin who helped me to
correct and improve the initial text of the manuscript.
This work was done in the stimulating atmosphere of
MSRI.

\end{document}